\documentclass[12pt,oneside]{amsart}
\usepackage{amsmath,amssymb} 
\newtheorem{thm}{Theorem}
\newtheorem{lem}[thm]{Lemma}

\theoremstyle{definition}
\newtheorem{defn}[thm]{Definition}

\newtheorem*{algorithm1}{Algorithm I}
\newtheorem*{algorithm2}{Algorithm II}
\newtheorem*{algorithm3}{Algorithm III}
\newtheorem*{algorithm4}{Algorithm IV}
\newtheorem*{algorithm5}{Algorithm V}

\newcommand{\bi}{\begin{itemize}}
\newcommand{\ei}{\end{itemize}}
\newcommand{\be}{\begin{enumerate}}
\newcommand{\ee}{\end{enumerate}}
\newcommand{\bc}{\begin{center}}
\newcommand{\ec}{\end{center}}
\newcommand{\bt}{\begin{tabular}}
\newcommand{\et}{\end{tabular}}

\newcommand{\abs}[1]{\left|#1\right|}

\newcommand{\hnn}{HNN-extension}
\newcommand{\Z}{\mathbb Z}

\newcommand{\ssm}{\smallsetminus}

\newcommand{\set}[1]{\left\{#1\right\}}

\def\Area{\hbox{\rm Area}}
\def\onto{{\kern3pt\to\kern-8pt\to\kern3pt}}

\def\ss{\smallskip}
\def\ms{\medskip}
\def\bs{\bigskip}
\def\ni{\noindent}

\begin{document}

\title[An isoperimetric function for Stallings' group]
{An isoperimetric function \\ for Stallings' group}

\author[W.~Dison]{W.Dison} \address{W.Dison, Department of Mathematics,
Imperial College London, London, SW7 2AZ, U.K.}
\email{\emph{william.dison@imperial.ac.uk}}

\author[M.~Elder]{M.Elder} \address{M.Elder, Department of Mathematical Sciences,
Stevens Institute of Technology, Hoboken, NJ 07030, USA} \email{\emph{melder@stevens.edu}}

\author[T.R.~Riley]{T.R.Riley} \address{T.R.Riley, Department of Mathematics, 310 Malott Hall, Cornell University, Ithaca, NY 14853-4201, USA}
\email{\emph{tim.riley@math.cornell.edu}}


\begin{abstract}
We prove that $n^{7/3}$ is an  isoperimetric function for a group of Stallings that is finitely presented but not of type $\mathcal{F}_3$.
\end{abstract}

\maketitle

\newlength{\miniwidth} \addtolength{\miniwidth}{\linewidth} \addtolength{\miniwidth}{-\parindent}

\section{Introduction}\label{sec:intro}

\ni
In the early 1960s Stallings \cite{Stallings} constructed a group $S$ enjoying the finiteness property $\mathcal{F}_2$ but not $\mathcal{F}_3$.  (A group is of type $\mathcal F_1$ when it can be finitely generated, $\mathcal F_2$ when it can be finitely presented, and more generally $\mathcal F_n$ when it admits an Eilenberg-Maclane space with finite $n$-skeleton.)  Bieri~\cite{Bieri} recognised $S$ to be
\begin{equation} \label{kernel}
\textup{Ker}( \, F(\alpha, \beta) \times F(\gamma, \delta) \times F(\epsilon, \zeta) \ \onto \ \Z \, )
\end{equation}
where the map is that from the product of three rank--$2$ free groups to $\Z = \langle t \rangle$ which sends all six generators to $t$, and he showed that using $(F_2)^n$ in place of $(F_2)^3$ gives a family of groups (the \emph{Bieri--Stallings} groups) of type $\mathcal{F}_{n-1}$ but not
$\mathcal{F}_n$ \cite{Bieri}.

Isoperimetric functions (defined below) for $S$ have been investigated by a number of authors.
Gersten proved that for $n \geq 3$, the groups in this family admit quintic
isoperimetric functions \cite{Gersten5}; this was sharpened to cubic by Baumslag, Bridson, Miller \& Short in the case of $S$ \cite[\S6]{BBMS}.
Bridson~\cite{Bridson4} argued that whenever $G_1$ and $G_2$ are finitely presentable groups admitting quadratic isoperimetric functions and epimorphisms $\phi_i:G_i \onto \Z$, if one doubles $G_1 \times G_2$ along the kernel of the map $\phi:  G_1 \times G_2 \onto \Z$, defined by $\phi(g_1,g_2) =\phi_1(g_1)+\phi_2(g_2)$, then the resulting group also admits a quadratic isoperimetric function. The Bieri--Stallings groups are examples of such \emph{doubles}.  But Groves found an error in his proof \cite{BridsonPersonal, GrovesPersonal}, and it seems that  Bridson's approach, in fact, gives cubic isoperimetric functions, generalising the result in \cite{BBMS}. In this article we prove:

\begin{thm} \label{Main Theorem}
Stallings' group $S$ has $n^{7/3}$ as an isoperimetric function.
\end{thm}

If the Dehn function of $S$ is not quadratic (i.e.\ not $\simeq n^2$, in the sense defined below) then it would be the first example of a subgroup of a $\textup{CAT}(0)$ group, namely $(F_2)^3$, with Dehn function bounded above by a polynomial, but not $\simeq n^{\alpha}$ for any $\alpha \in \mathbb{Z}$ --- we thank N.Brady for pointing this out.  Also, it would be an example of such a Dehn function `occurring naturally' rather than in a group especially constructed for the purpose such as in \cite{BB, SBR}.     If, on the other hand, the Dehn function of $S$ is quadratic then it would show the class of groups with quadratic Dehn functions to be wild enough to contain groups that are not of type $\mathcal{F}_3$, fulfilling Bridson's aim in \cite{Bridson4}.

Our theorem makes no reference to a specific finite presentation since, as is well--known, if such an isoperimetric function holds for one finite presentation of a group then it holds for all.  We will work with the presentation
\begin{equation} \label{pres}
 \langle \; a,b,c,d,s \; \mid \; [a,c], \; [a,d], \; [b,c], \; [b,d], \; s^a=s^b=s^c=s^d  \; \rangle
\end{equation}
for $S$ of \cite{BBMS, Gersten5}, with  $s^a=s^b=s^c=s^d$  shorthand for the six defining relations $s^as^{-b}$,  $s^as^{-c}$, $s^as^{-d}$, $s^bs^{-c}$,  $s^bs^{-d}$,  $s^cs^{-d}$.  One can view $S$ as an \hnn\ of the product of free groups $F(a,b) \times F(c,d)$ with stable letter $s$ commuting with all elements represented by words on $a^{\pm 1}, b^{\pm 1}, c^{\pm 1}, d^{\pm 1}$ of zero exponent-sum.
Gersten~\cite{Gersten5} shows this is related to the expression for  $S$ as a kernel (\ref{kernel}) via $a= \epsilon \alpha ^{-1},  b= \epsilon \beta ^{-1},  c =  \epsilon \gamma^{-1},  d= \epsilon \delta^{-1}, s= \zeta \epsilon^{-1}$.

Essentially, our strategy for establishing an $n^{7/3}$ isoperimetric function is to interplay two approaches to reducing words $w$ representing $1$ to $\varepsilon$. Both involve identifying a suitable subword $s^{\pm 1} \tau s^{\mp 1}$ of $w$, where $\tau=\tau(a,b,c,d)$, then converting $\tau$ to a word $\hat{\tau}$ in which the letters alternate between positive and negative exponent, and then cancelling off the $s^{\pm 1}$ with the $s^{\mp 1}$.  Repeating until all $s^{\pm 1}$ have been eliminated gives a word on $a,b,c,d$ that represents $1$ in $F(a,b) \times F(c,d)$.

In the first of these two approaches (Algorithm~I) all $a^{\pm 1},
b^{\pm 1}$ are shuffled to the start of $\tau$, leaving all the
$c^{\pm 1}, d^{\pm 1}$ at the end, and then letters $a^{\pm 1},
c^{\pm 1}$ are inserted to achieve the word $\hat{\tau}$ in the
required \emph{alternating} form.  The cost (see below) of
converting $\tau$ to $\hat{\tau}$ in this way is potentially great:
it can be as much as $\sim \! \ell(\tau)^2$; however control on the
length of $\hat{\tau}$ is good: it is always no more than $3
\ell(\tau)$.

The second approach (Algorithm~III) is to work through $\tau$ from
left to right inserting letters  $a^{\pm 1}, c^{\pm 1}$ as necessary
to achieve alternating form.  The cost of this algorithm and the
length of its output are heavily dependent on the internal structure
of $\tau$, but if $\tau$ possesses certain properties then good
bounds can be found.

In both cases the cost of cancelling off the $s^{\pm 1}, s^{\mp 1}$ is $\sim \! \ell(\hat{\tau})$.  Used alone, either approach would lead to a cubic  isoperimetric inequality.

\bs \ni \emph{Basic definitions}. $[x,y]:=x^{-1}y^{-1}xy$,  $x^y:= y^{-1}xy$, $x^{-y}:= y^{-1}x^{-1}y$. Write $u=u(a_1, \ldots, a_k)$ when $u$ is a word on the letters ${a_1}^{\pm 1}, \ldots, {a_k}^{\pm 1}$.
The length of $u$ as a word (with no free reductions performed) is $\ell(u)$.  The total number of occurrences of letters ${a_1}^{\pm 1}, \ldots, {a_l}^{\pm 1}$ in $u$ is $\ell_{a_1, \ldots ,a_l}(u)$.  Unless otherwise indicated, we consider two words to be the equal when they are identical letter-by-letter.

Given words $w, w'$ representing the same element of a group with finite presentation $\langle \mathcal{A} \mid \mathcal{R} \rangle$, one can convert $w$ to $w'$ via a sequence of words $W=(w_i)_{i=0}^m$ in which $w_0=w$, $w_m=w'$ and for each $i$, $w_{i+1}$ is obtained from $w_i$ by free reduction ($w_i = \alpha a a^{-1} \beta \mapsto \alpha \beta =w_{i+1}$ where $a \in \mathcal{A}^{\pm 1}$), by free expansion (the inverse of  a free reduction), or by applying a relator ($w_i=\alpha u\beta \mapsto \alpha v\beta = w_{i+1}$ where a cyclic conjugate of $uv^{-1}$  is in $\mathcal{R}^{\pm 1}$).   The \emph{cost} of $W$ is the number of $i$ such that
$w_i \mapsto w_{i+1}$ is an \emph{application-of-a-relator} move.  If  $w$ represents the identity (i.e.\ \emph{is null-homotopic}) then
 $\Area(w)$ is defined to be the minimal cost amongst all $W$ converting $w$ to the empty word $\varepsilon$, and the
\emph{Dehn function} $\Area : \mathbb{N} \to \mathbb{N}$ of $\langle \mathcal{A} \mid \mathcal{R} \rangle$ is
$$\Area(n) := \max \set{ \Area(w) \mid w = 1 \textup{ in } \Gamma \textup{ and } \ell(w) \leq n}.$$
An \emph{isoperimetric function} for $\langle \mathcal{A} \mid \mathcal{R} \rangle$ is any $f:  \mathbb{N} \to \mathbb{N}$  such that there exists $K>0$ for which $\Area(n) \leq K \, f(n$) for all $n$.
(The constant $K$ is not used by all authors, but is convenient for us here.)

For  $f,g : \mathbb{N} \to  \mathbb{N}$, we write $f \preceq g$ when $\exists C>0,  \forall n \in \mathbb{N}, f(n) \leq C g(Cn+C) + C n +C$, and we say $f \simeq g$ when $f \preceq g$ and $g \preceq f$.

\bs \ni \emph{Article organisation}. We give a number of
definitions, lemmas and algorithms in Section~\ref{prelims}.  In
Section~\ref{Skeletal proof} we use these to prove Theorem~\ref{Main
Theorem}.

\bs \ni \emph{Acknowledgements.}  We thank Noel~Brady,
Martin~Bridson, Daniel~Groves and especially Steve~Pride for many
fruitful discussions.  The third author is grateful for support from
NSF grant  DMS--0540830 and for the hospitality of the Institut des
Hautes \'Etudes Scientifique in Paris during the writing of this
article.

\section{Preliminaries.} \label{prelims}

\ni Our proof of Theorem~\ref{Main Theorem} will involve three classes of words.

\begin{defn} (\emph{Alternating} words.) \label{alternating}
A word $u = u(a,b,c,d)$ is \emph{alternating} if it is a concatenation of words $xy^{-1}$ in which $x,y \in \set{a,b,c,d}$.
\end{defn}

\ni [The reader familiar with van~Kampen diagrams and corridors (also known as bands) may find it helpful to note that alternating words are those which, after removing all $aa^{-1}, bb^{-1}, cc^{-1}$ and $dd^{-1}$ subwords, can be read along the sides of  $s$-corridors in van~Kampen diagrams over $S$.]

\begin{defn} (\emph{Balanced} words.)  \label{balanced}
A word $u=u(a,b,c,d,s)$ is \emph{balanced} if it has exponent sum zero and  in $S$ it represents an element of the subgroup $\langle a,b,c,d \rangle$.
\end{defn}

The following algorithm converts a word
$u=u(a,b,c,d)$ of exponent-sum zero into an alternating word of a \emph{preferred} form that represents the same element of $S$.

\begin{algorithm1}${}$ \newline
\indent \begin{minipage}{\miniwidth} Input a word of exponent-sum
zero $u=u(a,b,c,d)$.
 \begin{enumerate}
 \item Shuffle the letters $a^{\pm 1},b^{\pm 1}$ to the start of $u$ and freely reduce to give a word $\mu\lambda$ where $\mu= \mu(a,b)$ and  $\lambda = \lambda(c,d)$.
\item Intersperse $c^{\pm 1}$ through $\mu$ to give a word $\bar{\mu}=\bar{\mu}(a,b,c)$, and $a^{\pm 1}$ through $\lambda$ to give $\bar{\lambda} = \bar{\lambda}(a,c,d)$, such that:
\begin{enumerate}
\item $\bar{\mu}$ and $\bar{\lambda}$ are alternating,
\item for all $1 \leq i < \ell(\bar{\mu})/2$, exactly one of the $(2i-1)$-st and $(2i)$-th letters in $\bar{\mu}$ is $c$ or $c^{-1}$, and
\item for all $1 \leq j < \ell(\bar{\lambda})/2$, exactly one of the $(2j-1)$-st and $(2j)$-th letters in $\bar{\mu}$ is $a$ or $a^{- 1}$.
\end{enumerate}
\item Let $\kappa$ be the exponent-sum of $\mu$.   Insert $(ac^{-1})^{\kappa}$ between $\bar{\mu}$ and $\bar{\lambda}$.
\end{enumerate}
\ni Output $\bar{\mu} (ac^{-1})^{\kappa} \bar{\lambda}$.
\end{minipage}
\end{algorithm1}

\begin{defn} (\emph{Preferred} alternating words.)
A word $v=v(a,b,c,d)$ is in \emph{preferred alternating form} if it there is some $u$ such that the  output of Algorithm~I on input $u$ is $v$.
\end{defn}

\begin{lem} \label{length of paws}
The output $\bar{\mu} (ac^{-1})^{\kappa} \bar{\lambda}$ of Algorithm~I has length at most $3\ell(u)$ and $u$ can be converted to $\bar{\mu} (ac^{-1})^{\kappa} \bar{\lambda}$ at a cost of at most $10\ell(u)^2$.
\end{lem}

\begin{proof}
The exponent sum of $\mu$ is $\kappa$ and so that of $\lambda$ is $-\kappa$.  So $\abs{2 \kappa} \leq \ell(\mu \lambda) \leq \ell(u)$ and   $\ell(\bar{\mu} (ac^{-1})^{\kappa} \bar{\lambda}) \leq \ell(\bar{\mu} \bar{\lambda}) + \abs{2 \kappa} \leq 3 \ell(u)$.

The (crude) upper bound of $10 \ell(u)^2$ on the cost of converting $u$ to $\bar{\mu} (ac^{-1})^{\kappa} \bar{\lambda}$ holds because both $u$ and $\bar{\mu} (ac^{-1})^{\kappa} \bar{\lambda}$ can be converted to $\mu \lambda$ by shuffling letters and freely reducing at costs of at most $\ell(u)^2$ and
$\ell(\bar{\mu} (ac^{-1})^{\kappa} \bar{\lambda})^2 \leq (3\ell(u))^2$, respectively.
\end{proof}

The following lemma reveals balanced words to be those representing elements of the subgroup of $F(a,b) \times F(c,d)$ commuting with $s$ in the HNN-presentation of $S$.   We denote the centraliser of $s$ in $S$ by $C_S(s)$.

\begin{lem}
A word $u=u(a,b,c,d,s)$ represents an element $g$ of $$\langle a,b,c,d \rangle \ \cap \ C_S(s)$$ in $S$ if and only if $u$  is balanced.
\end{lem}

\begin{proof}
Note that all the relations of presentation (\ref{pres}) have
exponent-sum zero, so this quantity is preserved whenever a relation
is applied to a word.  Thus, if two words on the letters $a, b, c,
d, s$ represent the same element in $S$ then they have the same
exponent sum.

A word $u$ represents an element of $\langle a, b, c, d \rangle \cap
C_S(s)$ if and only if there exists an alternating word $v = v(a,b,
c, d)$ with $u$ and $v$ representing the same element of $S$. Since
any word on the letters $a, b, c, d$ with exponent-sum zero can be
converted into an alternating word by an application of Algorithm~I,
this is if and only if there exists a word $v = v(a, b, c, d)$ with
$u=v$ in $S$ and with $v$ having exponent-sum zero.  And by the
above remark this is if and only if $u$ represents an element of
$\langle a, b, c, d \rangle$ and itself has exponent-sum zero.
\end{proof}

\begin{lem} \label{balanced - alternating}
Suppose word $v_0 = v_0(a,b,c,d,s)$ is  expressed as $v_0=\alpha v_1 \beta$ in which $v_1$ is a balanced subword.  Then $v_0$  is balanced if and only if $\alpha\beta$ is balanced.
\end{lem}

\begin{proof}
Induct on $\ell_s(v_0)$, with the base case $\ell_s(v_0)=0$ immediate and the induction step  an application of Britton's Lemma.  Alternatively, this result is  an observation on the layout of $s$-corridors in a van~Kampen diagram demonstrating that $v_0$ equates to some alternating word in $S$.
\end{proof}

The next lemma concerns the existence of  balanced subwords within prescribed length-bounds in balanced words.

\begin{lem} \label{balanced existence}
If $\mu =\mu(a,b,c,d,s)$ is a balanced word with $\ell(\mu) \geq 4$, then for all $k \in [4, \ell(\mu)]$ there is a balanced subword $u$ of $\mu$ with $k/2 \ \leq \ \ell(u) \ \leq \ k$.
\end{lem}

\begin{proof}
We induct on $\ell(\mu)$.
First we identify certain balanced subwords $\alpha$ and $\beta$ in $\mu$.

\ss \ni \emph{Case: $\mu$ starts with a letter $x = s^{\pm 1}$.}  By Britton's Lemma, $\mu=x\alpha y \beta$ for $y=x^{-1}$ and for some balanced subword $\alpha$.

\ss \ni \emph{Case: $\mu$ starts with a letter $x \neq s^{\pm 1}$.}  Set a counter to $0$, then read through $\mu$ from left to right altering the counter as follows.  If the letter being read is not $s^{\pm 1}$ then add the exponent of that letter to the counter.  If it is $s^{\pm 1}$ then by Britton's Lemma that  $s^{\pm 1}$ is the first letter of a subword $s^{\pm 1}\gamma s^{\mp 1}$ such that $\gamma$ is balanced; hold the counter constant throughout $s^{\pm 1}\gamma s^{\mp 1}$ and then continue as before.  As $\mu$ is balanced, the counter will return to $0$ on reading some letter $y \neq s^{\pm 1}$ of opposite exponent to $x$.  Accordingly, $\mu=x \alpha y \beta$ in which $\alpha$ is balanced.

\ms

\ni In both cases, as $\alpha$ is balanced, so is $x \alpha y$, and hence so is $\beta$.

\ss
Now, in the base case of the induction we have $\ell(\mu)=4$ and so $k=4$, and we can take $u=\mu$. Indeed, whenever $\ell(\mu)=k$ we can take $u=\mu$, so let us assume henceforth that $\ell(\mu)>k$.

For the induction step, first suppose $\beta \neq \varepsilon$.  If $\max (\ell(x \alpha y), \ell(\beta)) \geq k$  then as $\ell(x \alpha y), \ell(\beta)<\ell(\mu)$ we can apply the induction hypothesis to obtain $u$.  If $\max (\ell(x \alpha y), \ell(\beta)) < k$ then both $x \alpha y$ and $\beta$ have length less than $k$ and, as $k < \ell(\mu)=\ell(x\alpha y)+\ell(\beta)$, either $\ell(x \alpha y)$ or $\ell(\beta)$ is at least $k/2$ and so serves as $u$.

Finally suppose $\beta = \varepsilon$.  If $\ell(\alpha) \geq k$ then, as $\ell(\alpha) < \ell(\mu)$, the induction hypothesis gives us $u$.  If $\ell(\alpha) < k$ then $\alpha$ serves as $u$ because  $\ell(\alpha) = \ell(\mu) -2 > k-2 \geq k/2$ since $k \geq 4$.
\end{proof}

The remainder of this section works towards Algorithm~IV which will convert a balanced word $u=u(a,b,c,d,s)$ into a preferred alternating word $v$ representing the same element of $S$. 

The next algorithm concerns converting a word $\tau_0 = \tau_0(a,b,c,d)$ into alternating form by working through it from left to right inserting letters $a^{\pm 1}$ as needed.  In contrast to Algorithm~III, which is an elaboration of this algorithm, the group element represented will not be preserved.  The purpose of this algorithm is to define a number $P(\tau_0)$ which will be the difference in length of the input and output words, a quantity which plays an important role in our analysis of Algorithm~III.

\begin{algorithm2} ${}$ \newline
\indent \begin{minipage}{\miniwidth}
Input a word $\tau = \tau(a,b,c,d)$.   Define $\beta_i$ to be the length-$(\ell(\tau)-i)$ suffix of $\tau$. Define $\tau_0:= \tau$ and $\alpha_0 := \varepsilon$.

The algorithm will produce a sequence of words $(\tau_i)_{i=0}^{\ell(\tau)}$ of the form $\tau_i = \alpha_i \beta_i$,  where $\alpha_i$ is an alternating word or an alternating word  concatenated with an $a,b,c,d$.   For $0 \leq i <  \ell(\tau)$, obtain $\alpha_{i+1}$ from $\alpha_i$ as follows.  We have $\beta_i = x \beta_{i+1}$ for some letter $x$.

\ms \ni Case: $\ell(\alpha_i)$ is even.
\begin{itemize}
\item If $x \in \set{a,b,c,d}$, then  $\alpha_{i+1} :=  \alpha_i x$.
\item If $x \in \set{a^{-1},b^{-1}, c^{-1},d^{-1}}$, then $\alpha_{i+1} :=  \alpha_i a x$.
\end{itemize}

\ms \ni Case: $\ell(\alpha_i)$ is odd.
\begin{itemize}
\item If $x \in \set{a,b,c,d}$, then  $\alpha_{i+1} := \alpha_i  a^{-1} x$.
\item If $x \in \set{a^{-1},b^{-1},c^{-1},d^{-1}}$, then $\alpha_{i+1} = \alpha_i x$.
\end{itemize}

\ms
\ni Output $\tau_{\ell(\tau)}$.
\end{minipage}
\end{algorithm2}

\begin{defn}
For words $\tau = \tau(a,b,c,d)$, define $P(\tau)$ to be the number of letters $a^{\pm 1}$ inserted by Algorithm~II on input $\tau$.
\end{defn}

\begin{lem} \label{P with alt subwords removed}
Let $\Pi$ be a collection of $p$ disjoint alternating subwords of a
word $\tau = \tau(a, b, c, d)$, and let $\bar{\tau}$ be the word
formed from $\tau$ by removing all the subwords specified by $\Pi$.
Then $P(\tau) \leq P(\bar{\tau}) + 2p$.
\end{lem}

\begin{proof}
For a letter $l \in \{a^{\pm1}, b^{\pm1}, c^{\pm1}, d^{\pm1} \}$
write $\chi(l) \in \{\pm1\}$ for the exponent of $l$. For a word $w
= w(a, b, c, d)$ write $w[i]$ for the $i^\text{th}$ letter of $w$.
For $i \in \{2, \ldots, \ell(w)\}$ define $$d_i(w) \ = \
\begin{cases}
1 \quad \text{if $\chi(w[i]) = \chi(w[i-1])$,}\\
0 \quad \text{if $\chi(w[i]) \not= \chi(w[i-1])$,}
\end{cases}$$  and define $$d_1(w) \ = \ \begin{cases}
  1 \quad \text{if $\chi(w[1]) = -1$,}\\
  0 \quad \text{if $\chi(w[1]) = 1$.}
\end{cases}$$  Note that, during the running of Algorithm~II on a word
$\tau$, an $a^{\pm1}$ is inserted during the transition from
$\tau_{i-1}$ to $\tau_i$ precisely when $d_i(w) = 1$.  Thus $P(\tau)
= \sum_{i=1}^{\ell(\tau)} d_i(\tau)$.

By induction, it suffices to prove the lemma in the case $p=1$.
Suppose $\tau = uvw$ and $\bar{\tau} =uw$ for some words $u, v, w$
with $v$ alternating. Note that: \begin{align*}
  &d_i(\tau) = d_i(u) = d_i(\bar{\tau}) &&i = 1, \ldots, \ell(u);\\
  &d_{\ell(u)+i}(\tau) = d_i(v) = 0 && i = 2, \ldots, \ell(v);\\
  &d_{\ell(uv)+i}(\tau) = d_i(w) = d_{\ell(u)+i}(\bar{\tau}) &&i =
  2, \ldots, \ell(w).
\end{align*}  Thus \begin{align*}
  P(\tau) \ &= \ \sum_{i=1}^{\ell(\tau)} d_i(\tau)\\
  &= \ \left[\sum_{i=1}^{\ell(u)} d_i(\bar{\tau}) \right] + d_{\ell(u)+1}(\tau) +
  d_{\ell(uv)+1}(\tau) + \left[ \sum_{i=2}^{\ell(w)}
  d_{\ell(u)+i}(\bar{\tau}) \right]\\
  &\leq \ \left[ \sum_{i=1}^{\ell(\bar{\tau})}
  d_i(\bar{\tau})\right]
  +d_{\ell(u)+1}(\tau) + d_{\ell(uv) + 1}(\tau)\\
  &\leq \ P(\bar{\tau}) +2.
\end{align*}
\end{proof}

\begin{defn}
For words $\sigma = \sigma(a,b,c,d,s)$, define $Q(\sigma)$ to be the number of times letters $b^{\pm 1}$ alternate with letters $d^{\pm 1}$ in $\sigma$.  More precisely, if  the word obtained from $\sigma$ by deleting all letters $a^{\pm 1}, c^{\pm 1}, s^{\pm 1}$ is $\mu_1 \nu_1\mu_2 \nu_2\ldots \mu_q \nu_q$, in which $\nu_1,\mu_2,\nu_2,\ldots, \nu_{q-1}, \mu_q \neq \varepsilon$ and $\mu_i = \mu_i(b)$ and $\nu_i = \nu_i(d)$ for all $i$,  then $Q(\sigma) = q$.  
\end{defn}

\begin{defn}
For words $\sigma = \sigma(a,b,c,d,s)$, define $R(\sigma)$ to be the maximum over all suffixes $\beta$ of $\sigma$ of the absolute value of the exponent sum of $\beta$.
\end{defn}

The following algorithm works through a word $\tau_0 = \tau_0(a,b,c,d)$ from left to right inserting letters  $a^{\pm1}, c^{\pm1}$ without changing the element of $F(a,b) \times F(c,d)$ it represents.  If the exponent-sum of $\tau$ is zero then the output will be alternating.

\begin{algorithm3}${}$ \newline
\indent \begin{minipage}{\miniwidth}
Input a word $\tau_0 = \tau_0(a,b,c,d)$.   Define $\beta_i$ to be the length-$(\ell(\tau)-i)$ suffix of $\tau_0$. Define $\alpha_0 := \varepsilon$ and $\Delta_0 := \varepsilon$.

The algorithm will produce a sequence of words $(\tau_i)_{i=0}^{\ell(\tau)}$ of the form $\tau_i = \alpha_i \Delta_i \beta_i$,  where $\alpha_i$ is an alternating word or an alternating word  concatenated with an $a,b,c,d$, and $\Delta_i$ is $a^r$ or $c^r$ for some $r  \in \mathbb{Z}$.

For $0 \leq i <  \ell(\tau)$, obtain $\alpha_{i+1}, \Delta_{i+1}$ from $\alpha_i, \Delta_i$ as follows.  We have $\beta_i = x \beta_{i+1}$ for some letter $x$.

\ms \ni Case (1) \   $\Delta_i = a^r$ for some $r \in \mathbb{Z} \ssm \set{0} $  and $\ell(\alpha_i)$ is even.
\begin{itemize}
\item[(1.1)] If $x \in \set{a,c,d}$ then  $\alpha_{i+1} :=  \alpha_i x$ and $\Delta_{i+1} := \Delta_i$.
\item[(1.2)] If $x = b$ then  $\alpha_{i+1} = \alpha_i (ac^{-1})^r  x$ and  $ \Delta_{i+1} := c^r$.
\item[(1.3)] If $x \in \set{a^{-1},c^{-1},d^{-1}}$ then $\alpha_{i+1} :=  \alpha_i a x$ and $\Delta_{i+1} :=  a^{r-1}$.
\item[(1.4)] If $x = b^{-1}$ then  $\alpha_{i+1} =  \alpha_i  (ac^{-1})^r  c x$ and $\Delta_{i+1} := c^{r-1}$.
\end{itemize}

\ms \ni Case (2) \   $\Delta_i =   a^r$  for some $r \in \mathbb{Z}  \ssm \set{0}$ and $\ell(\alpha)$ is odd.
\begin{itemize}
\item[(2.1)] If $x \in \set{a,c,d}$ then  $\alpha_{i+1} := \alpha_i  a^{-1} x$ and $\Delta_{i+1} := a^{r+1}$.
\item[(2.2)] If $x = b$ then  $\alpha_{i+1} := \alpha_i (c^{-1}a)^r c^{-1} x$ and $\Delta_{i+1} :=  c^{r+1}$.
\item[(2.3)] If $x \in \set{a^{-1},c^{-1},d^{-1}}$ then $\alpha_{i+1} = \alpha_i x$ and  $\Delta_{i+1} := \Delta_i$.
\item[(2.4)] If $x = b^{-1}$ then $\alpha_{i+1} := \alpha_i  (c^{-1}a)^r  x$ and $\Delta_{i+1} :=  c^r$.
\end{itemize}

\ms \ni When  $\Delta_i =  c^r$  for some $r \in \mathbb{Z}$, obtain $\alpha_{i+1}$ and $\Delta_{i+1}$  similarly, but with $a,b$ interchanging roles with $c,d$.  Call the cases involved (3.1--3.4) and (4.1--4.4).

\ms \ni Output $\tau_{\ell(\tau)}$.
\end{minipage}
\end{algorithm3}

\begin{lem} \label{to alternating}
Suppose $\tau = \tau(a,b,c,d)$ is  a word of exponent sum zero.  Then Algorithm~III converts $\tau$ to an alternating word $\hat{\tau}$ with $\ell(\hat{\tau}) \geq \ell(\tau)$, with $Q(\hat{\tau}) = Q(\tau)$, and with
\begin{equation} \label{length control}
\ell(\hat{\tau}) - \ell(\tau) \ \leq \   P(\tau) + 4 (R(\tau) +1) Q(\tau).
\end{equation}
Moreover, the cost of transforming $\tau$ to $\hat{\tau}$ is at most
\begin{equation} \label{cost}
(R(\tau) +2)\ell(\tau) + 2 (R(\tau) +2)^2 Q(\tau).
\end{equation}
\end{lem}

\begin{proof}
As the exponent sum of each $\tau_{i+1}$ is the same as that of $\tau_i$, it remains at zero throughout the run of the algorithm and $\Delta_{\ell(\tau)}$ must be $\varepsilon$.  It follows that $\hat{\tau} = \tau_{\ell(\tau)}$ is alternating.

If one removes all letters $a^{\pm 1}$ and $c^{\pm 1}$ from  $\hat{\tau}$ and $\tau$, they become identical words, and so   $Q(\hat{\tau}) = Q(\tau)$.

The transformation $\tau_i$ to  $\tau_{i+1}$ can be achieved  at a cost of at most $\abs{r}+1$ in Cases~$\star$.1, $\star$.3 and at most  $(\abs{r}+1)^2$ in Cases~$\star$.2, $\star$.4, where $r$ is the exponent in $\Delta_i$. In every instance,
\begin{equation} \label{exp sum}
\abs{r} \ \leq \ R(\tau) + 1.
\end{equation}

Suppose removing all letters $a^{\pm 1}$ and $c^{\pm 1}$ from $\tau$
gives $\mu_1 \nu_1\mu_2 \nu_2\ldots \mu_q \nu_q$, in which
$\nu_1,\mu_2,\nu_2,\ldots, \nu_{q-1}, \mu_q \neq \varepsilon$ and
$\mu_i = \mu_i(b)$ and $\nu_i = \nu_i(d)$ for all $i$.  By
definition, $Q(\tau)=q$.  The process described above will carry a
power of $c$ through the word from the left until it hits $\nu_1$,
when it will be converted to a power of $a$, which will then be
carried until it hits $\mu_2$ when it reverts to a power of $c$, and
so on.   So Cases~$\star$.2 and $\star$.4 are invoked either $2
Q(\tau)-1$ or $2 Q(\tau) - 2$ times depending on whether or not
$\nu_q  = \varepsilon$. Cases~$\star$.1 and $\star$.3 are invoked
the remaining  $\ell(\tau) - 2  Q(\tau) + 1$ or $\ell(\tau) - 2
Q(\tau) + 2$ times.  Combining these estimates we see that the total
cost of converting $\tau$ to $\hat{\tau}$ is at most
\begin{equation}
(R(\tau)+2) (\ell(\tau) - 2 Q(\tau)+2)  +  ( R(\tau) + 2)^2 (2 Q(\tau)-1),
\end{equation}
which, discarding some negative terms and noting that $Q(\tau) \geq
1$, gives (\ref{cost}).

The length estimate (\ref{length control}) comes from counting the
letters deposited into $\tau$ in the above process  en route to
reaching $\hat{\tau}$.  They occur in two forms.   (Note: we do not
consider the powers of  $a$ and $c$ carried through the word as
\emph{deposited}.)  Firstly, there are the single letters $a^{\pm
1}$ or  $c^{\pm 1}$ inserted in Cases 1.3,  1.4,  2.1, 2.2, 3.3,
3.4, 4.1, and 4.2.  These total $P(\tau)$.  And, secondly, there are
the $(ac^{-1})^{r}$ of Cases~1.2 and 1.4, the $(c^{-1}a)^{r}$ of 2.2
and 2.4, the $(ca^{-1})^{r}$ of 3.2 and 3.4, and the $(a^{-1}c)^{r}$
of 4.2 and 4.4. These cases occur less than $2Q(\tau)$ times and by
(\ref{exp sum}) each inserts a word of length at most $2\abs{r} \leq
2(R(\tau) + 1)$.
\end{proof}


Our next algorithm transforms a balanced word on $a,b,c,d,s$ to a word in preferred alternating form that represents the same element of $S$.

\begin{algorithm4} ${}$ \newline
\indent \begin{minipage}{\miniwidth}
Input a balanced word $u =u(a,b,c,d,s)$. Define $u_0:= u$ and $L := \ell_s(u)/2$.  Then for $0\leq i < L$ recursively obtain $u_{i+1}$ from $u_i$ by the following two steps.

\begin{itemize}
\item[(A)]  Locate a subword $s^{\pm 1} \tau_i s^{\mp 1}$ in $u_i$  such that  $\tau_i = \tau_i(a,b,c,d)$ and has zero exponent sum (which, by  Britton's Lemma, we know exists).  Use Algorithm~III to transform $\tau_i$ into alternating form $\hat{\tau}_i$.

\item[(B)] Shuffle the $s^{\pm 1}$ through  $\hat{\tau}_i$ and cancel it with the $s^{\mp 1}$ to give  $u_{i+1}$.
\end{itemize}

\ni This produces $u_L$, which contains no letters $s^{\pm 1}$.  Next --

\begin{itemize}

\item[(C)] Reverse every instance of Step~A  to get a word $\bar{u}$, which is $u$ with all letters $s^{\pm 1}$ deleted.

\item[(D)] Run Algorithm~I on $\bar{u}$ to give a word $v$ in preferred alternating form.
\end{itemize}
\ni Output $v$.
\end{minipage}
\end{algorithm4}

\begin{lem} \label{counting the cost}
Suppose $u =u(a,b,c,d,s)$ is a balanced word with  $2 \leq \ell(u)
\leq m$.  Suppose $\Pi$ is some collection of at most $p$ disjoint
subwords in $u$, each in preferred alternating form, such that
deleting these subwords leaves a word of length at most $k$.  Then
Algorithm~IV transforms $u$ into $v$  at a cost of no more than
$$80k^3 + 75k^2p + 16m^2.$$
\end{lem}

\begin{proof}
By Lemma~\ref{to alternating}, performing Step~A on $\tau_i$ costs at most
\begin{equation} \label{cost first est}
 (R(\tau_i)   +   2)\ell(\tau_i) + 2 (R(\tau_i) +2)^2 Q(\tau_i).
\end{equation}

Now, for all $i$,
\begin{eqnarray}
 Q(\tau_i) & \leq & Q(u_i) \  = \  Q(u)  \ \leq \ k+p, \label{bound on Q}
\end{eqnarray}
as letters $b^{\pm 1}$ alternate with letters $d^{\pm 1}$ at most once in each preferred--alternating--form subword of $u$, and transformations as per Lemma~\ref{to alternating} do not alter $Q$.

Note that if $ \sigma'$ is obtained from a word $\sigma$ by deleting
a collection of disjoint alternating subwords, then $R(\sigma) \leq
R( \sigma') +1$.  We recursively define a collection $\Pi_i$ of
disjoint alternating subwords of $u_i$ in the letters $a^{\pm1},
b^{\pm1}, c^{\pm1}, d^{\pm1}$ by $\Pi_0:= \Pi$ and for $0 \leq i
<L$, $$\Pi_{i+1} := (\Pi_i \ssm \Pi'_i) \cup \set{ \hat{\tau}_i}.$$
Let $\Pi'_i$ by the subset of $\Pi_i$ consisting of those words
which have letters in common with $\tau_i$.  Note that each word in
$\Pi'_i$ is a subword of $\tau_i$ since it contains no occurrence of
a letter $s^{\pm1}$.  Removing the subwords $\Pi'_i$ from $\tau_i$
produces a word $\tau'_i$ whose letters all originate in  $u$ but
not in any of its subwords $\Pi$.  If  $i \neq j$ then $\tau'_i$ and
$\tau'_j$ originate from different letters in $u$ so
\begin{eqnarray} \label{sum length tau'}
  \sum_{i=0}^{L-1} \ell(\tau'_i) \ \leq \ k.
\end{eqnarray} Furthermore, since $L \leq k/2$ and $ R(\tau'_i)\leq \ell(\tau'_i)$
one has
\begin{eqnarray}\label{Rsum}
  \sum_{i=0}^{L-1} R(\tau_i)  \  \leq \ \sum_{i=0}^{L-1}
  (R(\tau'_i) + 1) \ \leq \ \frac{k}{2} + \sum_{i=0}^{L-1}
  \ell(\tau'_i)    \  \leq \ \frac{3k}{2}.
\end{eqnarray}

Similarly \begin{eqnarray} \label{P of tau}
  P(\tau_i) \ \leq \ \ell(\tau'_i) + 2|\Pi'_i|
  \ \leq \ \ell(\tau'_i) + 2|\Pi_i| \ \leq \ \ell(\tau'_i) + 2p + k
\end{eqnarray} for all $i$, by Lemma~\ref{P with alt subwords removed}
applied to removing the subwords $\Pi'_i$ from $\tau_i$ and noting
that $|\Pi_{j+1}| \leq |\Pi_j| + 1$ for all $j$ and that $i \leq L
\leq k/2$. It follows from \ref{P of tau} and \ref{sum length tau'}
that
\begin{eqnarray} \label{sum of P less than k}
  \sum_{i=0}^{L-1} P(\tau_i) \ \leq \ k + pk + \frac{k^2}{2}
  \ \leq \ 2k^2 + pk.
\end{eqnarray}

Now, for all $j$,
\begin{eqnarray}
\ell(u_{j+1}) - \ell(u_j) & = & \ell(\hat{\tau}_j) -  \ell( \tau_j) \nonumber \\
& \leq & P(\tau_j) + 4(R(\tau_j)+1)Q(\tau_j) \nonumber \\
& \leq & P(\tau_j) + 4(R(\tau_j)+1)(k+p), \label{telescope}
\end{eqnarray}
where the first and second inequalities are applications of
(\ref{length control}) and (\ref{bound on Q}), respectively.  So,
for all $i \leq L$,
\begin{eqnarray}
  \ell(u_i) & \leq & m + \sum_{j=0}^{i-1}
  \left( P(\tau_j) + 4(R(\tau_j)+1)(k+p) \right) \nonumber \\
  & \leq & m + 2k^2 + pk + (6k + 4i)(k+p) \nonumber \\
  & \leq & m + 2k^2 + pk + 8k(k + p) \nonumber \\
  & \leq & m + 10k^2 + 9kp \label{length of ui}
\end{eqnarray}
where the first inequality uses $\ell(u_0) = m$ and
(\ref{telescope}), the second uses (\ref{Rsum}) and (\ref{sum of P
less than k}), and the third uses the inequality $i \leq L \leq
k/2$.  We thus find
\begin{eqnarray} \label{sum of lengths}
  \sum_{i=0}^{L-1} \ell(u_i) \ \leq \ \frac{1}{2}km + 5k^3 + \frac{9}{2}k^2p.
\end{eqnarray}

The total cost of all instances of Step~A, including those
implemented within Step~C, is at most
\begin{align}
  2\sum_{i=0}^{L-1} & \left[ (R(\tau_i) + 2)\ell(\tau_i)
  + 2 (R(\tau_i) +2)^2 Q(\tau_i)\right]  \nonumber \\
  & \leq \ 2 \max_{i} (\ell(\tau_i)) \sum_{i=0}^{L-1} (R(\tau_i) + 2)
  + 4 \max_{i} Q(\tau_i) \left[ \sum_{i=0}^{L-1} (R(\tau_i) + 2) \right]^2  \nonumber \\
  & \leq \ 2(m + 10k^2 + 9kp) \left( \frac{3}{2} k + k \right)
  + 4(k+p)\left( \frac{3}{2} k + k \right)^2 \nonumber \\
  & \leq \ 75k^3 + 70k^2p + 5mk
\end{align}
where the initial estimate comes from (\ref{cost first est}), and
the second inequality uses $ L \leq k/2$,  $\ell(\tau_i)  \leq
\ell(u_i)$, (\ref{bound on Q}),  (\ref{Rsum}), and (\ref{length of
ui}).

The total cost of all instances of Step~B is $\sum_{i=0}^{L-1}
\ell(\hat{\tau}_i)$ which, by (\ref{sum of lengths}), is at most
$\frac{1}{2}km + 5k^3 + \frac{9}{2}k^2p$, as $\ell(\hat{\tau}_i)
\leq \ell(u_i)$. As $\ell(\bar{u}) \leq m$, Lemma~\ref{length of
paws} tells us that the cost of Step~D is at most $10m^2$. Summing
these three cost estimates gives the total cost as at most
\begin{align}
  (\frac{1}{2}km &+ 5k^3 +\frac{9}{2}k^2p) + (75k^3 + 70k^2p + 5mk) + 10m^2 \nonumber \\
  & \quad \leq 80k^3 + \frac{149}{2}k^2p + \frac{11}{2}km + 10m^2 \nonumber \\
  & \quad \leq 80k^3 + 75k^2p + 16m^2.
\end{align} where for the final inequality we have used that $k \leq
m$.
\end{proof}

\section{Proof of Theorem~\ref{Main Theorem}.}  \label{Skeletal proof}

\ni Our final algorithm concerns converts null-homotopic words in
$S$ to $\varepsilon$.  Its cost analysis will establish
Theorem~\ref{Main Theorem}.  (The finitely many $w$ of length less
than $8$ are irrelevant for the asymptotics of the Dehn function of
$S$.)  It constructs a sequence of null-homotopic words
$(w_i)_{i=0}^l$ and a subset $T_i$ of $\set{1,2, \ldots, \ell(w_i)}$
specifying a collection of the letters of $w_i$ by their locations.

\begin{algorithm5}  ${}$ \newline
\indent \begin{minipage}{\miniwidth} Input a word $w$ of length at
least $8$ representing $1$ in $S$. Let  $n := \ell(w)$. Input a
parameter $k \in [4, n]$.

Define $w_0:=w$ and let $T_0$ be the empty set.  For successive $i$
such that $T_i \neq \set{1,2, \ldots, \ell(w_i)}$, obtain $w_{i+1}$
and $T_{i+1}$ from $w_i$ and $T_i$  by performing the following
steps.

\begin{itemize}
\item[(i)] Let $\bar{w}_i$ be the word obtained from $w_i$ by deleting the letters
in the positions $T_i$.    If  $\ell(\bar{w}_i) \leq k$, then define
$u := w_i$.  If $\ell(\bar{w}_i) > k$ then let $\bar{u}$ be a
balanced subword of $\bar{w}_i$ with $k/2 \leq \ell(\bar{u}) \leq k$
and take $u$ to be the longest subword of $w_i$ which reduces to
$\bar{u}$ when all the letters specified by $T_i$ are removed.  In
either case, $u$ is a balanced subword of $w_i$ of which between
$k/2$ and $k$ letters are not in positions in $T_i$.

\item[(ii)]
Use Algorithm~IV to convert $u$ to a word $v$ in preferred alternating form.   Obtain $w_{i+1}$ from $w_i$ by replacing $u$ with $v$.  Let  $T_{i+1}$ be the locations in $w_{i+1}$ of the letters of $v$ and  of the letters that originate in $w_i$ (but not in $u$) and have locations in $T_i$.
\end{itemize}

\ni Define $l := i+1$ where $i$ is the value in the final run of Steps (i) and (ii).

\begin{itemize}
\item[(iii)]
Reduce $w_l$ to $\varepsilon$ by shuffling the letters $a^{\pm 1}, b^{\pm 1}$ to the start of the word and then freely reducing.
\end{itemize}
\end{minipage}
\end{algorithm5}

\ni \emph{Notes.}  Steps (i) and (ii) are repeated at most $((n-k)/(k/2))+1 = (2n/k)-1$ times and so  $l \leq 2n/k$.

The reason $w_i = w_{i+1}$ for all $0 \leq i < l$ is that the words $u$ and $v$ in Step~(ii) represent the same element of $S$.  One effect of each instance of Step~(ii) is to remove any letters $s^{\pm 1}$ in $u$.   In particular none of the letters of $w_i$ specified by $T_i$ are $s^{\pm 1}$ and so there are no $s^{\pm 1}$ in $w_l$.

For all $i$,  the letters of $w_i$ specified by  $T_i$ comprises $\leq i$ subwords, each in preferred alternating form. (The  number of these subwords rises by at most one with each run of Steps (i) and (ii).)

The existence of the $u$ of Step~(i) follows from
Lemmas~\ref{balanced - alternating} and \ref{balanced existence} as
follows.   As $w_i$ is null-homotopic and hence balanced,
Lemma~\ref{balanced - alternating} applies and tells us that
$\bar{w}_i$ is also balanced, since $T_i$ specifies a number of
alternating subwords in $w_i$.  As $\ell \left(\bar{w}_i\right) > k
\geq 4$, Lemma~\ref {balanced existence} applies and tells us that
$\bar{w}_i$ contains a balanced subword $\bar{u}$ with $k/2 \leq
\ell(\bar{u}) \leq k$.  An appeal to Lemma~\ref{balanced -
alternating} tells us that $u$ is balanced. In the final case, where
$u$ is $w_{l-1}$, we see that $u$ is balanced because it is
null-homotopic.

The viability of Step (iii) follows from the facts that  $w_l$ is null-homotopic in $S$ and contains no letters $s^{\pm 1}$, and so is null-homotopic in $$F(a,b) \times F(c,d) \ = \ \langle \; a,b,c,d \; \mid \; [a,c], \; [a,d], \; [b,c], \; [b,d] \; \rangle.$$

\ms \ni \emph{Cost analysis.} For all $i$, the letters of $w_i$
specified by $T_i$ comprise a number of subwords in preferred
alternating form.  The number of such subwords specified by $T_i$ is
at most $i \leq l \leq 2n/k$, since the number increases by at most
one for each transition $T_j$ to $T_{j+1}$. Each specified subword
has length at most three times the length of a corresponding subword
of $w$, by  Lemma~\ref{length of paws}. Thus $\ell(w_i) \leq 3
\ell(w)$ and so, for each $i$, the subword $u$ of step (i) has
$\ell(u) \leq 3\ell(w)$.  Applying Lemma~\ref{counting the cost}
with $p = 2 n / k$ and $m=3n$ tells us that the cost of each
instance of Step~(ii) is at most
\begin{equation*}
  80k^3 + 150nk + 144n^2.
\end{equation*}
Multiplying  by $l \leq 2n/k$ and adding $(3n)^2$, which is an upper
bound on the cost of Step (iii) as $\ell(w_l) \leq 3n$, gives the
estimate $$\Area(w) \ \leq \ 160nk^2 + 309n^2 + 288\frac{n^3}{k}.$$
So $k = n^{2/3}$, which is compatible with the condition $k \geq 4$
since we assumed that $n \geq 8$, gives our $n^{7/3}$ isoperimetric
function.   \qed

\bibliography{bibli}
\bibliographystyle{plain}

\end{document}